# Mathematical Analysis of the Role of Information on the Dynamics of Typhoid Fever

**Nyanga H. Masasila[1,2], Rigobert C. Ngeleja[3], Odeli J. Kigodi[4]**

[1]Department of General Studies, Dar es Salaam Institute of Technology, Dar es Salaam, Tanzania
[2]Doctoral School of Informatics, University of Debrecen, Debrecen, Hungary
[3]Commission for Science and Technology (COSTECH), Dar es Salaam, Tanzania
[4]Department of Mathematics, Physics and Informatics, Mkwawa University College of Education, University of Dar es Salaam, Iringa, Tanzania
Email: hondanyanga@gmail.com

## Abstract

This study presents a deterministic model to examine how information affects the spread of Typhoid Fever. The model's properties, including its stability and basic reproduction number, are analyzed. Simulations show that information can influence behavior in ways that may increase disease transmission. Notably, the rise in Typhoid cases is linked to poor adherence to health precautions. The findings highlight the critical role of public education in controlling the disease and emphasize the need to include information campaigns in prevention strategies.

## 1. Introduction

Typhoid fever is an exclusively human enterically transmitted systemic disease caused by infection with the bacterium Salmonella enterica serovar Typhi. Although largely controlled in Europe and North America, typhoid remains endemic in many parts of the world, notably Africa, where it is an important cause of febrile illness in crowded, low-income settings [1]. The infection is often passed through contaminated food and drinking water, and it is more prevalent in the places where handwashing is less frequent. Moreover, it can be transmitted to the susceptible human being through adequate contact with the infected person.

This infection produces bacteraemic illness, with prolonged high fever, headache and malaise being characteristic symptoms [2]. Other symptoms might include confusion, diarrhea and vomiting. Without effective treatment, typhoid fever can lead to the altered mental states and fatal at large [2] [3]. The only treatment for typhoid is antibiotics (the commonly used are ciprofloxacin and ceftiaxone). To prevent its transmission the policies address that before traveling to the high-risk areas, vaccination against typhoid fever is mandatory.

Typhoid is among the most endemic diseases, and thus of major public health concern in tropical developing counties like Tanzania [3]. Therefore, the **information** about its spread and transmission in a given area becomes stimulant of awareness and prevention among the people. The government and other authorities use social media and related means to circulate the information on eruption of Typhoid.

Several studies have applied mathematical modeling and optimal control techniques to understand disease dynamics and economic systems. For instance, Uwakwe *et al.* [4] developed an optimal control model to study the transmission dynamics of Avian Spirochaetosis, incorporating time-dependent controls to minimize infection and cost. Similarly, Katende [5] applied mathematical modeling to analyze the growth dynamics of infant financial markets, highlighting the effectiveness of nonlinear systems in describing real-world phenomena. While these studies demonstrate the utility of advanced modeling tools, they do not consider the influence of information dissemination as a behavioral control mechanism, particularly in the context of typhoid fever.

Despite extensive research on the mathematical modeling of infectious diseases, most existing models of typhoid fever focus predominantly on medical interventions such as treatment, vaccination, and environmental sanitation. These models often neglect the role of public awareness and behavioral responses driven by information dissemination, an increasingly important factor in disease control, especially in the digital age. Moreover, while optimal control methods have been applied to various epidemiological and economic models, their application in modeling information-sensitive control strategies for typhoid fever remains largely unexplored. This gap is significant because behavioral change influenced by timely and accurate information can alter the course of an epidemic without solely relying on biomedical interventions. Therefore, this study addresses a critical void by integrating information-driven behavioral dynamics with optimal control theory, providing a more comprehensive and realistic framework for understanding and managing typhoid fever outbreaks. In this paper, therefore, we present the mathematical model which explains the effect of information on dynamics of Typhoid fever. Furthermore, we introduce essential parameters that can lead to the reduction of the spread of diseases based on the information received by susceptible people.

# 2. Model development

## 2.1. Model Description

This Typhoid Model is in two settings, the human beings and the transmitting bacteria in the environment that is referred to food and water denoted by $B$. We divide the Human population into three subgroups: first is a group of people who have not acquired the infection but may get it if they adequately get into contact with infectious human $I$ or infectious media (environment) $A$ to be known as susceptible and symbolized by $S$, second are the infectious human being who can transmit the disease symbolized by $I$, when individual from subgroup $I$ get treated or through strong body immunity may recover and attain a temporary immunity known as recovered population symbolized by $R$ otherwise they die naturally at a rate $\pi_2$ or because of the disease at the rate $\pi_3$. The transmitting bacteria in the environment which includes objects, food or water contaminated with *Salmonella enterica serotype Typhi bacteria* also play as an agent of transmission of typhoid fever if they get into adequate contact with a susceptible human being.

## 2.2. Description of Interaction

When the susceptible Human being come into contact with the infectious agent the dynamics begins. A human being may be infected through eating or drinking contaminated food or water that has pathogens at the rates $\theta_1$ (fecal-oral transmission). Moreover, human beings may be infected through adequate contact with other infectious human being at a rate $\theta_2$. Human beings are recruited at a constant rate $\pi_1$ and removed by naturally death at the rate $\pi_2$. If not treated human beings may die due to the disease at a rate $\pi_3$. The bacteria causing typhoid that are in the environment (food or water) are recruited constantly at a rate $\eta_1$ and removed at a rate $\lambda_3$. Additionally, Salmonella enterica serotype Typhi bacteria may also be populated to the environment by the infected human beings ($I$) at the rate $\eta_2$.

## 2.3. Variable and Parameters and Their Description (Table 1)

Table 1. Parameters and their description for Typhoid fever.

| Parameters | Description | Value | Source |
|---|---|---|---|
| $\lambda$ | Immunity loss rate of $R$ | 0.1255 | Estimated |
| $\pi$ | Recruitment rate of human population | 0.92 | Estimated |
| $\theta$ | Adequate contact rate: $S$ and $A$ | 0.95 | [6] |
| $\eta$ | Shading rate of bacteria | 0.95 | [6] |
| $\rho$ | Information induced behaviour response of $S$ | 0.07 | Estimated |
| $C$ | Concentration of bacteria in the environment | 000 | [7] |
| $\gamma(B)$ | Probability of a human being catching Dysentery | 0.0001 | Estimated |
| $\pi$ | Death rate of human beings | 0.005 | Estimated |

**Continued**

| | | | |
|---|---|---|---|
| $\lambda$ | Recovery rate | 0.048 | Estimated |
| $\pi$ | Disease-induced death rate for human beings | 0.015 | [8] |
| $\lambda$ | Removal of rate of $A$ | 0.025k | Estimated |
| $\nu$ | Rate information spread which depends on $I$ | 0.025k | Estimated |
| $\theta$ | Adequate contact rate: $I$ and $S$ | 0.0021 | [9] |
| $\eta$ | Recruitment of bacilli in $A$ | 0.95 | Estimated |

where; $\gamma(B)=\frac{B}{B+C}$.

## 3. Model Assumption

The typhoid disease model is developed based on the assumption below:

1) Human population who are susceptible are recruited at a constant rate;

2) The natural death rate of all human beings in this model is the same;

3) Human population mix homogeneously;

4) Human beings from all subgroups have equal chance of being infected by Typhoid fever;

5) All other media that can transfer the disease are included in one compartment called $B$.

Considering the dynamics illustrated in the model development and the stated model assumption, we can summarize the dynamics of typhoid fever in a compartmental diagram in Figure 1. It captures the interaction between human beings and Salmonella enterica serotype Typhi bacteria in the environment (food and water).

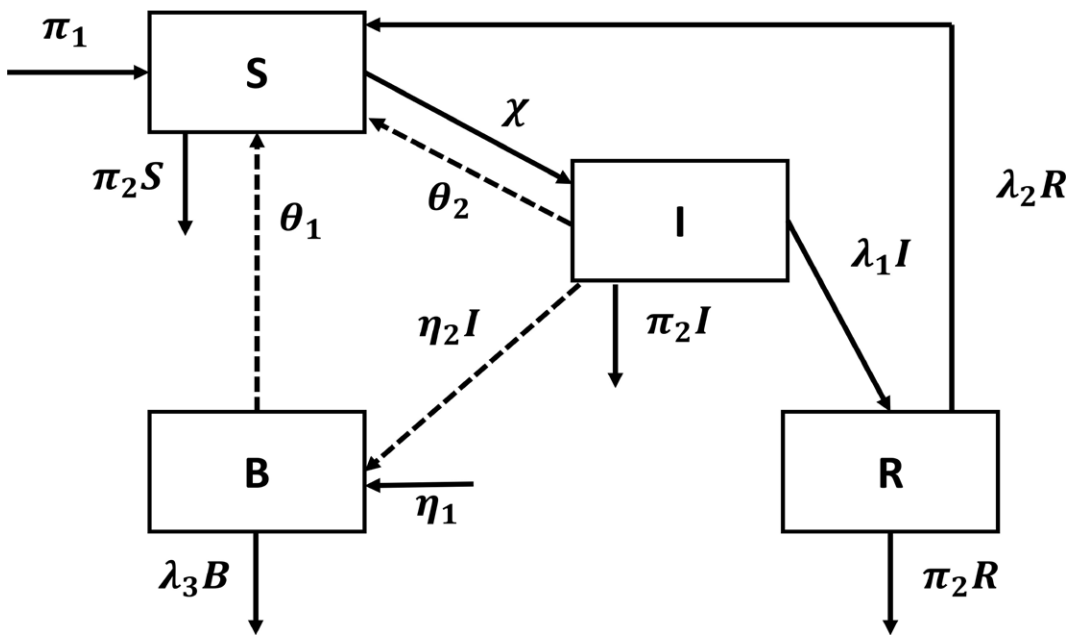

Figure 1. Compartmental model for Typhoid Fever.

where $\chi=\theta_1(1-\rho)\gamma(B)S(t)+\frac{\theta_2(1-\rho)I(t)S(t)}{1+\nu_b I}$

### Model Equation

The dynamics of Typhoid fever is represented by system (1)

$$\frac{\mathrm{d}S(t)}{\mathrm{d}t}=\pi_1+\lambda_2 R(t)-\theta_1(1-\rho)\gamma(B)S(t)-\frac{\theta_2(1-\rho)I(t)S(t)}{1+\nu_b I}-\pi_2 S(t), \quad \text{(1a)}$$

$$\frac{\mathrm{d}I(t)}{\mathrm{d}t}=\theta_1(1-\rho)\gamma(B)S(t)+\frac{\theta_2(1-\rho)I(t)S(t)}{1+v_bI}-(\lambda_1+\pi_2+\pi_3)I(t), \quad (1b)$$

$$\frac{\mathrm{d}R}{\mathrm{d}t}=\lambda_1I(t)-(\pi_2+\lambda_2)R(t), \quad (1c)$$

$$\frac{\mathrm{d}B}{\mathrm{d}t}=\eta_1+\eta_2\frac{I(t)}{N}-\lambda_3B(t). \quad (1d)$$

$$N=S+I+R$$

# 4. Properties of the Typhoid Fever Model

## 4.1. Model's Invariant Region

In the modelling process of Typhoid in human beings we assume that the model's state variables and parameters are non-negative for $\forall t\geq 0$. The model system is then analyzed in an appropriate, feasible region that satisfies this assumption. Using Theorem 4.1 stated below, we can obtain the invariant region of the typhoid fever model as given below.

**Theorem 4.1** *All forward model solutions in* $R_+^4$ *of the typhoid fever model system are feasible* $\forall t\geq 0$ *if they go in the invariant region* $\Lambda$ *for* $\Lambda=\omega_1\times\omega_2$ given that;

$$\omega_1=(S,I,R)\in R_+^3:S+I+R\leq N$$
$$\omega_2=B\in R_+^1$$

thus, the positive invariant region of Typhoid fever system is symbolized as $\Lambda$.

*Proof.* **For human population:**

Here, we need to show that the solutions of the typhoid fever model system (1) are feasible $\forall t>0$ as they enter region $\omega_1$ for human population.

Let $\omega_1=(S,I,R)\in R^3$ be solution space of the typhoid fever model system with positive initial conditions.

We will then have,

$$\frac{\mathrm{d}N}{\mathrm{d}t}=\frac{\mathrm{d}S}{\mathrm{d}t}+\frac{\mathrm{d}I}{\mathrm{d}t}+\frac{\mathrm{d}R}{\mathrm{d}t} \quad (2)$$

Substituting the system equations into (2) yields,

$$\frac{\mathrm{d}N}{\mathrm{d}t}=\pi_1-\pi_2N-\pi_3I$$

It then gives

$$\frac{\mathrm{d}N}{\mathrm{d}t}\leq\pi_1-\pi_2N$$

Which then yields,

$$\frac{\mathrm{d}N}{\mathrm{d}t}+\pi_2N\leq\pi_1$$

Using the integrating factor method,

we use $IF=\mathrm{e}^{\pi_2t}$ which when multiplied throughout gives

$$\mathrm{e}^{\pi_2 t}\frac{\mathrm{d}N}{\mathrm{d}t}+\mathrm{e}^{\pi_2 t}N\pi_2 \le \pi_1 \mathrm{e}^{\pi_2 t}$$

which gives

$$\frac{\mathrm{d}\left(N\mathrm{e}^{\pi_2 t}\right)}{\mathrm{d}t} \le \pi_1 \mathrm{e}^{\pi_2 t}$$

Integrating on both sides yields

$$N\mathrm{e}^{\pi_2 t} \le \frac{\pi_1}{\pi_2}\mathrm{e}^{\pi_2 t}+C$$

This then gives;

$$N \le \frac{\pi_1}{\pi_2}+C\mathrm{e}^{-\pi_2 t}$$

We then plug in $t=0, N(t=0)=N_0$ as the initial conditions which yield;

$$N_0-\frac{\pi_1}{\pi_2} \le C$$

Then, the substitution of the constant gives,

$$N \le \frac{\pi_1}{\pi_2}+\left(N_0-\frac{\pi_1}{\pi_2}\right)\mathrm{e}^{-\pi_2 t}$$

When $N_0 > \frac{\pi_1}{\pi_2}$, human population are asymptotically reduced to $\frac{\pi_1}{\pi_2}$ and when $N_0 < \frac{\pi_1}{\pi_2}$ the human population are asymptotically enlarged to $\frac{\pi_1}{\pi_2}$.

This then proves that all feasible solutions of the typhoid fever model system for human population go into the region

$$\omega_1 = \left\{(S,I,R): N \le \max\left\{N_0, \frac{\pi_1}{\pi_2}\right\}\right\}$$

**For Bacteria in the media(environment)**

In this section, we also need to show that the solutions of the typhoid fever system for the bacteria in the media are feasible $\forall t>0$ whenever they go into an invariant region $\omega_2$. With non-negative initial condition, we now let the solution of the system to be $\omega_2 = B \in R_+^1$

from the equation

$$\frac{\mathrm{d}B}{\mathrm{d}t} = \eta_1+\eta_2\frac{I(t)}{N}-\lambda_3 B(t). \tag{3}$$

But

$$I \le N$$

Then this implies that

$$\frac{I}{N} \le 1.$$

Substituting into Equation (3) we obtain;

$$\frac{dB}{dt} \le \eta_1 + \eta_2 - \lambda_3 B.$$

It then gives

$$\frac{dB}{dt} + \lambda_3 B \le \eta_1 + \eta_2.$$

Using the integrating factor method we will have

$$IF = e^{\lambda t}$$

Then

$$e^{\lambda_3 t}\frac{dB}{dt} + e^{\lambda_3 t}\lambda_3 B \le e^{\lambda_3 t}(\eta_1 + \eta_2).$$

$$\frac{d\left(Be^{\lambda_3 t}\right)}{dt} \le (\eta_1 + \eta_2)e^{\lambda_3 t}.$$

$$Be^{\lambda_3 t} \le \frac{\eta_1 + \eta_2}{\lambda_3}e^{\lambda_3 t} + C,$$

$$B(t) \le \frac{\eta_1 + \eta_2}{\lambda_3} + Ce^{-\lambda_3 t}.$$

We then use $t = 0, B(t = 0) = B_0$ as the initial conditions which gives

$$B_0 - \frac{\eta_1 + \eta_2}{\lambda_3} \le C,$$

$$B(t) \le \frac{\eta_1 + \eta_2}{\lambda_3} + \left(B_0 - \frac{\eta_1 + \eta_2}{\lambda_3}\right)e^{-\lambda_3 t}.$$

When $B_0 > \frac{\eta_1 + \eta_2}{\lambda_3}$ the concentration of bacteria in the environment are asymptotically reduced to $\frac{\eta_1 + \eta_2}{\lambda_3}$ and when $B < \frac{\eta_1 + \eta_2}{\lambda_3}$ the concentration of bacteria in the environment asymptotically enlarged to $\frac{\eta_1 + \eta_2}{\lambda_3}$.

This then proves that all feasible solutions of the typhoid fever model system for bacteria in the environment go into the region

$$\Omega_B = \left\{B : B \le \max\left\{B_0, \frac{\eta_1 + \eta_2}{\lambda_3}\right\}\right\}$$

## 4.2. Positivity of the Solution

In this section, we are required to show that the variables and parameters used in the typhoid fever model are greater than or equal to zero $\forall t \ge 0$.

**Theorem 4.2** *We assume the initial values of the system* (1) *to be:* $S(0) > 0$ *and* $(I(0), R(0), B(0)) \ge 0$. *Then the solution set of the typhoid fever model system such that* $S(t), I(t), R(t)$ *and* $B(t)$ *are positive* $\forall t \ge 0$.

*Proof.* Here, the requirements are to show that the solution of each individual equation form of the Typhoid system (1) is positive

Now consider Equation (1a) of the typhoid fever system,

$$\frac{\mathrm{d}S(t)}{\mathrm{d}t}=\pi_1+\lambda_2 R(t)-\theta_1(1-\rho)\gamma(A)S(t)-\frac{\theta_2(1-\rho)I(t)S(t)}{1+\nu_b I}-\pi_2 S(t)$$
$$\geq -\left(\theta_1(1-\rho)\gamma(A)S+\frac{\theta_2(1-\rho)I(t)}{1+\nu_b I}+\pi_2\right)S$$
$$\frac{\mathrm{d}S}{\mathrm{d}t}\geq -\left(\theta_1(1-\rho)\gamma(A)S+\frac{\theta_2(1-\rho)I(t)}{1+\nu_b I}+\pi_2\right)S$$

Integration yields

$$S\geq S_0\mathrm{e}^{-\int_0^t\left(\theta_1(1-\rho)\gamma(A)S+\frac{\theta_2(1-\rho)I(t)}{1+\nu_b I}+\pi_2\right)\mathrm{d}\tau}>0$$

since

$$\mathrm{e}^{-\int_0^t\left(\theta_1(1-\rho)\gamma(A)S+\frac{\theta_2(1-\rho)I(t)}{1+\nu_b I}+\pi_2\right)\mathrm{d}\tau}>0.$$

Considering the next equation, we have;

$$\frac{\mathrm{d}I(t)}{\mathrm{d}t}=\theta_1(1-\rho)\gamma(B)S(t)+\frac{\theta_2(1-\rho)I(t)S(t)}{1+\nu_b I}-(\lambda_1+\pi_2+\pi_3)I(t)$$

Thus

$$\frac{\mathrm{d}I}{\mathrm{d}t}\geq -(\lambda_1+\pi_2+\pi_3)I.$$

Integration yields

$$I\geq I_0\mathrm{e}^{-(\lambda_1+\pi_2+\pi_3)t}>0$$

since

$$\mathrm{e}^{-(\lambda_1+\pi_2+\pi_3)}>0.$$

Considering another equation of system (1) we will have

$$\frac{\mathrm{d}R}{\mathrm{d}t}=\lambda_1 I(t)-(\pi_2+\lambda_2)R(t)$$

Thus

$$\frac{\mathrm{d}R}{\mathrm{d}t}\geq -(\pi_2+\lambda_2)R.$$

Integrating we get

$$R\geq R_0\mathrm{e}^{-(\pi_2+\lambda_2)t}>0,$$

since

$$\mathrm{e}^{-(\pi_2+\lambda_2)}>0.$$

**For bacteria in the environment**

Here, we consider the last equation of the system (1) which is given below;

$$\frac{\mathrm{d}B}{\mathrm{d}t}=\eta_1+\eta_2\frac{I(t)}{N}-\lambda_3 B(t)$$

Then, we will have

$$\frac{\mathrm{d}B}{\mathrm{d}t}\geq -\lambda_3 B.$$

Integrating we get

$$B \geq B_0 e^{-\lambda_3 t} > 0,$$

Since

$$e^{-\lambda_3 t} > 0.$$

□

## 5. Analysis of the Model

Here, we work on the presence and stability of the stationary points and the conditions for extinction or persistence of the disease (basic reproduction number).

### 5.1. Disease-Free Equilibrium

In order to get the disease Free Equilibrium point, we set the variables $I$, $R$ and $B$ of the typhoid fever system equals zero, such that $I = R = 0$ and $B = 0$.

Now substituting the above into the system (1), we obtain the disease-free-equilibrium point of the typhoid system as given in (4)

$$X_0\left(S^0, I^0, R^0, B^0\right) = \left(\frac{\pi_1}{\pi_2}, 0, 0, 0\right). \tag{4}$$

### 5.2. Computation of the Basic Reproduction Number $R_0$

To find the basic reproduction number, we use the next-generation method by [10] [11]. Consider a heterogeneous population in compartments $S, I, R$ and $B$ arranged such that $m$ infectious classes come first. Assume $F_i(x)$ as rate of entrance of new individual with typhoid fever in class $i$, $V_i^+(x)$ rate of transfer of individuals in the class $i$ by any other means except the typhoid fever induced $V_i^-(x)$ be the rate of transfer of individuals out of class $i$. The basic reproduction number ( $R_0$ ) represents the average number of secondary typhoid fever cases caused by one infectious individual during their infectious period in a fully susceptible population.

- If $R_0 < 1$, each infected person transmits the disease to less than one person on average, suggesting that typhoid fever can die out, and the disease-free equilibrium is stable.
- If $R_0 > 1$, an infected person transmits the disease to more than one person, allowing the disease to persist in the population. This makes the disease-free equilibrium unstable.
- If $R_0 = 1$, each infected person infects exactly one other person, meaning the disease can remain in the population without causing a major outbreak (Allen *et al.* [12]).

To calculate $R_0$, the next-generation method (as in [10] [11]) is used. In a compartmental model with compartments $S$, $I$, $R$, and $B$, assume $m$ infectious classes come first. Let $F_i(x)$ represent the rate of new infections in class $i$, $V_i^+(x)$ the rate of entry into class $i$ from other causes, and $V_i^-(x)$ the rate of leaving class $i$.

The model system is as presented below;

$$x_i' = F_i(x) - V_i(x) \tag{5}$$

where $V_i(x) = V_i^-(x) - V_i^+(x)$.

Then, we use $x_0$, to find the $m \times m$ matrices $\boldsymbol{F}$ and $\boldsymbol{V}$

$$\boldsymbol{F} = \left(\frac{\partial F_i}{\partial x_j}(x_0)\right), \quad \boldsymbol{V} = \left(\frac{\partial V_i}{\partial x_j}(x_0)\right) \tag{6}$$

with $1 \le i, j \le m$.

By using the study by [13], we call Matrix $\boldsymbol{F}\boldsymbol{V}^{-1}$, a next-generation matrix and $\rho\left(\boldsymbol{F}\boldsymbol{V}^{-1}\right)$ is the basic reproduction number. Arranging the typhoid system starting with the infectious classes we get the rearranged system (7)

$$\frac{\mathrm{d}I(t)}{\mathrm{d}t} = \theta_1(1-\rho)\gamma(B)S(t) + \frac{\theta_2(1-\rho)I(t)S(t)}{1+\nu_b I} - (\lambda_1 + \pi_2 + \pi_3)I(t) \tag{7a}$$

$$\frac{\mathrm{d}B}{\mathrm{d}t} = \eta_1 + \eta_2\frac{I(t)}{N} - \lambda_3 B(t). \tag{7b}$$

$$\frac{\mathrm{d}R}{\mathrm{d}t} = \lambda_1 I(t) - (\pi_2 + \lambda_2)R(t), \tag{7c}$$

$$\frac{\mathrm{d}S(t)}{\mathrm{d}t} = \pi_1 + \lambda_2 R(t) - \theta_1(1-\rho)\gamma(A)S(t) - \frac{\theta_2(1-\rho)I(t)S(t)}{1+\nu_b I} - \pi_2 S(t) \tag{7d}$$

We consider the infectious classes (7a) to (7b) with compartment $I$ and $B$ from the system.

$$\boldsymbol{F}_i = \begin{pmatrix} \left(\theta_1(1-\rho)\gamma(B) + \dfrac{\theta_2(1-\rho)I(t)}{1+\nu_b I}\right)S \\ \eta_1 \end{pmatrix} \tag{8}$$

And

$$\boldsymbol{V}_i = \begin{pmatrix} (\lambda_1 + \pi_2 + \pi_3)I \\ \lambda_3 B - \eta_2\dfrac{I}{N} \end{pmatrix}. \tag{9}$$

We then get matrices of $\boldsymbol{F}$ and $\boldsymbol{V}$ which are the Jacobian matrices at $x_0$

$$\frac{\partial F_i}{\partial x_j} = \begin{pmatrix} \dfrac{\partial F_1}{\partial I} & \dfrac{\partial F_1}{\partial B} \\ \dfrac{\partial F_2}{\partial I} & \dfrac{\partial F_2}{\partial B} \end{pmatrix} = \begin{pmatrix} \left(\dfrac{\theta_2(1-\rho)}{(1+\nu_b I)^2}\right)S & \left(\dfrac{\theta_1(1-\rho)C}{(B+C)^2}\right)S \\ 0 & 0 \end{pmatrix}.$$

Now at $x_0$ we will have

$$\boldsymbol{F} = \begin{pmatrix} \dfrac{\theta_2(1-\rho)\pi_1}{\pi_2} & \dfrac{\theta_1(1-\rho)\pi_1}{C\pi_2} \\ 0 & 0 \end{pmatrix}. \tag{10}$$

Then we have;

$$\boldsymbol{V} = \frac{\partial V_i}{\partial x_j}(x_0) = \begin{pmatrix} \dfrac{\partial V_1}{\partial I} & \dfrac{\partial V_1}{\partial B} \\ \dfrac{\partial V_2}{\partial I} & \dfrac{\partial V_2}{\partial B} \end{pmatrix}$$

$$V=\begin{pmatrix}\lambda_1+\pi_2+\pi_3 & 0\\ -\frac{\eta_2\pi_2}{\pi_1} & \lambda_3\end{pmatrix}. \tag{11}$$

From (11) we can easily obtain $V^{-1}$ and $FV^{-1}$.

We then use maple to compute the basic reproduction number, which is given in (12);

$$R_0=\frac{(1-\rho)(C\pi_1\lambda_3\theta_2+\pi_2\eta_2\theta_1)}{\pi_2(\lambda_1+\pi_2+\pi_3)C\lambda_3} \tag{12}$$

## 6. Steady State and Stability of the Critical Points

We prove the presence and stability of the stationary points of the system (1).

### 6.1. Disease-Free Equilibrium

The disease-free-equilibrium point of the Typhoid fever model is as given in (13)

$$X_0\left(S^0,I^0,R^0,B^0\right)=\left(\frac{\pi_1}{\pi_2},0,0,0\right). \tag{13}$$

### 6.2. Local Stability of the Disease-Free Equilibrium Point

This section presents the analysis of local stability of the disease-free stationary point of the typhoid fever model. We use Jacobian method by considering that all equations in typhoid fever model in (1) are analyzed at the disease-free stationary point $X_0$. We are required to compute and assess the eigenvalues of Jacobian matrix ($J(X_0)$) to verify that the disease-free stationary point is locally and asymptotically stable. Furthermore, we need to show that the real parts of the eigenvalues of the matrix at $X_0$ are negative.

Using the concept by Martcheva [14], we are required to show that eigenvalues are negative, in which we need to prove that determinant of the Jacobian matrix is positive and its trace negative.

The matrix $J(X_0)$ at $X_0$ is given by:

$$J(X_0)=\begin{pmatrix}-\pi_2 & -\frac{\theta_2(1-\rho)\pi_1}{\pi_2} & \lambda_2 & -\frac{\theta_1(1-\rho)\pi_1}{C\pi_2}\\ 0 & -(\lambda_1+\pi_1+\pi_3) & 0 & \frac{\theta_1(1-\rho)\pi_1}{C\pi_2}\\ 0 & \lambda_1 & -(\pi_2+\lambda_2) & 0\\ 0 & \frac{\eta_2\pi_2}{\pi_1} & 0 & -\lambda_3\end{pmatrix} \tag{14}$$

The computation clearly shows that the trace of the matrix (14) is negative and given by

$$-(\pi_2+\lambda_1+\pi_1+\pi_3+\pi_2+\lambda_2+\lambda_3)$$

For the determinant of matrix (14), using maple software we can find the determinant of a Jacobian matrix as in (15):

$$\frac{\pi_2(\pi_2+\lambda_2)(\lambda_3 C\lambda_1+\lambda_3\pi_1 C+\lambda_3 C\lambda_2-\eta_2\theta_1(1-\rho))}{C} \tag{15}$$

which is positive if and only if $R_0 < 1$.

The above results justify that the typhoid free stationary point $X^0$ is locally asymptotically stable as in the theorem below:

**Theorem 6.1** *The Disease Free Equilibrium* $X_0$ *of Typhoid Fever is locally asymptotically stable if* $R_0 < 1$ *and unstable if* $R_0 > 1$.

### 6.3. Global Stability of the Disease-Free Equilibrium Point

Here we analyse the global stability of the disease-free equilibrium point. We use Metzler matrix method as stated by (Castillo-Chavez *et al.* *[15]*). To do this, we first subdivide the general system (1) of typhoid fever into transmitting and non-transmitting components.

Now let $\boldsymbol{Y}_n$ be the vector for non-transmitting compartment, $\boldsymbol{Y}_i$ be the vector for transmitting compartment and $\boldsymbol{Y}_{X_0,n}$ be the vector of disease-free point. Then

$$\begin{cases} \dfrac{\mathrm{d}\boldsymbol{Y}_n}{\mathrm{d}t} = \boldsymbol{A}_1\left(\boldsymbol{Y}_n - \boldsymbol{Y}_{X_0,n}\right) + \boldsymbol{A}_2\boldsymbol{Y}_i \\ \dfrac{\mathrm{d}\boldsymbol{Y}_i}{\mathrm{d}t} = \boldsymbol{A}_3\boldsymbol{Y}_i \end{cases} \tag{16}$$

We then have

$$\boldsymbol{Y}_n = (S,R)^{\mathrm{T}} \quad \boldsymbol{Y}_i = (I,B) \quad \boldsymbol{Y}_{X_0,n} = \left(\frac{\vartheta}{\mu}, 0\right)$$

$$\boldsymbol{Y}_n - \boldsymbol{Y}_{X_0,n} = \begin{pmatrix} S - \dfrac{\pi_1}{\pi_2} \\ R \end{pmatrix}$$

In order to prove that the Desease-free equilibrium point is globally stable we need to show that Matrix $\boldsymbol{A}_1$ has real negative eigenvalues and $\boldsymbol{A}_3$ is a Metzler matrix in which all off-diagonal elements must be non-negative. Referring to (16), we write the general model as below

$$\begin{pmatrix} \pi_1 + \lambda_2 R - \theta_1(1-\rho)\gamma(B)S - \dfrac{\theta_2(1-\rho)IS}{1+\nu_b I} - \pi_2 S \\ \lambda_1 I(t) - (\pi_2+\lambda_2)R(t). \end{pmatrix} = \boldsymbol{A}_1\begin{pmatrix} S - \dfrac{\pi_1}{\pi_2} \\ R \end{pmatrix} + \boldsymbol{A}_2\begin{pmatrix} I \\ B \end{pmatrix}$$

and

$$\begin{pmatrix} \theta_1(1-\rho)\gamma(B)S(t) + \dfrac{\theta_2(1-\rho)I(t)S(t)}{1+\nu_b I} - (\lambda_1+\pi_2+\pi_3)I(t) \\ \eta_1 + \eta_2\dfrac{I(t)}{N} - \lambda_3 B(t) \end{pmatrix} = \boldsymbol{A}_3\begin{pmatrix} I \\ B \end{pmatrix}$$

We then use the transmitting and non-transmitting elements from the Typhoid fever model to get the following matrices.

$$\boldsymbol{A}_1 = \begin{pmatrix} -\pi_2 & \lambda_2 \\ 0 & -(\pi_2+\lambda_2) \end{pmatrix} \tag{17}$$

$$\boldsymbol{A}_2 = \begin{pmatrix} -\dfrac{\theta_2(1-\rho)\pi_1}{\pi_2} & -\dfrac{\theta_1(1-\rho)\pi_1}{C\pi_2} \\ \lambda_1 & 0 \end{pmatrix} \tag{18}$$

$$\boldsymbol{A}_3 = \begin{pmatrix} \dfrac{\theta_2(1-\rho)\pi_1}{\pi_2} - (\lambda_1+\pi_2+\pi_3) & \dfrac{\theta_1(1-\rho)\pi_1}{C\pi_2} \\ \dfrac{\eta_2\pi_2}{\pi_1} & -\dfrac{\lambda_3}{C} \end{pmatrix} \tag{19}$$

Considering matrix $\boldsymbol{A}_1$, it is clear through computation that the eigenvalues are real and negative, which now confirms that the system

$$\frac{\mathrm{d}\boldsymbol{Y}_n}{\mathrm{d}t} = \boldsymbol{A}_1\left(\boldsymbol{Y}_n - \boldsymbol{Y}_{X_0,n}\right) + \boldsymbol{A}_2\boldsymbol{Y}_i$$

is globally and asymptotically stable at $\boldsymbol{Y}_{X_0}$.

Considering matrix $\boldsymbol{A}_3$ it is clear that all its off-diagonal elements are non-negative and thus $\boldsymbol{A}_3$ is a Metzler stable matrix. Therefore, Disease-Free Equilibrium point for Typhoid Fever system is globally asymptotically stable and as a result we have the following theorem:

**Theorem 6.2** *The disease-free equilibrium point is globally asymptotically stable in* $E_0$ *if* $R_0 < 1$ *and unstable if* $R_0 > 1$.

## 6.4. Existence of Endemic Equilibrium

In this section, we investigate the conditions that support the existence of the endemic equilibrium point of the system (1).

The endemic equilibrium point $E^*\left(S^*, I^*, R^*, B^*\right)$ is obtained by solving the equations obtained by setting the derivatives of (1) equal to zero. We then have system (20) which exists for $R_0 > 1$.

$$\pi_1 + \lambda_2 R(t) - \theta_1(1-\rho)\gamma(B)S(t) - \frac{\theta_2(1-\rho)I(t)S(t)}{1+v_b I} - \pi_2 S(t) = 0 \tag{20a}$$

$$\theta_1(1-\rho)\gamma(B)S(t) + \frac{\theta_2(1-\rho)I(t)S(t)}{1+v_b I} - (\lambda_1+\pi_2+\pi_3)I(t) = 0 \tag{20b}$$

$$\lambda_1 I(t) - (\pi_2+\lambda_2)R(t) = 0 \tag{20c}$$

$$\eta_1 + \eta_2\frac{I(t)}{N} - \lambda_3 B(t) = 0 \tag{20d}$$

We will prove its existence in the endemic equilibrium points of the Typhoid Fever using the approach described in the studies by [16] and [17]. For the endemic equilibrium to exist it must satisfy the condition $I \neq 0$ or $R \neq 0$ or $B \neq 0$ that is $S > 0$ or $R > 0$ or $I > 0$ or $B > 0$ must be satisfied. Now adding system (20), we have

$$\pi_1 - \pi_2(S+I+R) - \pi_3 I + \eta_1 + \eta_2\frac{I}{N} - \lambda_3 B(t) = 0 \tag{21}$$

But from Equation (20d), we have $\eta_1 + \eta_2 \dfrac{I}{N} - \lambda_3 B(t) = 0$

and $S + I + R = N$ It follows that

$$\pi_1 = \pi_2 N + \pi_3 I$$

Now since $\pi_1 > 0$, $\pi_2 > 0$ and $\pi_3 > 0$ we can discern that $\pi_2 N > 0$ and $\pi_3 I > 0$ implying that $S > 0$, $I > 0$, $R > 0$ and $M > 0$.

This proves that the endemic equilibrium point of Typhoid Fever disease exists.

### 6.5. Global Stability of Endemic Equilibrium Point

In this section, we determine the conditions under which the endemic equilibrium points are stable or unstable. In which we prove whether the solution starting sufficiently close to the equilibrium remains close to the equilibrium and approaches the equilibrium as $t \to \infty$, or if there are solutions starting arbitrarily close to the equilibrium which does not approach it respectively.

As postulated in the study by [10], we assert that the local stability of the Disease-Free Equilibrium advocates for local stability of the Endemic Equilibrium for the reverse condition. We thus find the global stability of Endemic equilibrium using a Korobeinikov approach as described by [10] [18] [19].

We formulate a suitable Lyapunov function for Typhoid Fever model as given in the form below:

$$V = \sum a_i \left( y_i - y_i^* \ln y_i \right)$$

where $a_i$ is defined as a properly selected positive constant, $y_i$ defines the population of the $i^{th}$ compartment, and $y_i^*$ is the equilibrium point.

We will then have

$$V = W_1 \left( S - S^* \ln S \right) + W_2 \left( I - I^* \ln I \right) + W_3 \left( R - R^* \ln R \right) + W_4 \left( B - B^* \ln B \right)$$

The constants $W_i$ are non-negative in $\Lambda$ such that $W_i > 0$ for $i = 1,2,3,4$. The Lyapunov function $V$ together with its constants $W_1, W_2, W_3, W_4$ chosen in such a way that $V$ is continuous and differentiable in a space.

We then compute the time derivative of $V$ from which we get:

$$\frac{\mathrm{d}V}{\mathrm{d}t} = W_1 \left( 1 - \frac{S^*}{S} \right) \frac{\mathrm{d}S}{\mathrm{d}t} + W_2 \left( 1 - \frac{I^*}{I} \right) \frac{\mathrm{d}I}{\mathrm{d}t} + W_3 \left( 1 - \frac{R^*}{R} \right) \frac{\mathrm{d}R}{\mathrm{d}t} + W_4 \left( 1 - \frac{B^*}{B} \right) \frac{\mathrm{d}B}{\mathrm{d}t}$$

Now using the Typhoid Fever system (1), we will have

$$\begin{aligned} \frac{\mathrm{d}V}{\mathrm{d}t} = {} & W_1 \left( 1 - \frac{S^*}{S} \right) \left[ \pi_1 + \lambda_2 R(t) - \theta_1 (1-\rho) \gamma(B) S(t) - \frac{\theta_2 (1-\rho) I(t) S(t)}{1 + \nu_b I} - \pi_2 S(t) \right] \\ & + W_2 \left( 1 - \frac{I^*}{I} \right) \left[ \theta_1 (1-\rho) \gamma(B) S(t) + \frac{\theta_2 (1-\rho) I(t) S(t)}{1 + \nu_b I} - (\lambda_1 + \pi_2 + \pi_3) I(t) \right] \end{aligned}$$

$$+W_3\left(1-\frac{R^*}{R}\right)\left[\lambda_1 I(t)-(\pi_2+\lambda_2)R(t)\right]$$
$$+W_4\left(1-\frac{I^*}{I}\right)\left[\eta_1+\eta_2\frac{I(t)}{N}-\lambda_3 B(t)\right]$$

At endemic equilibrium point after the substitution and simplification into time derivative of $V$, we get:

$$\frac{\mathrm{d}V}{\mathrm{d}t}=-W_1\left(1-\frac{S^*}{S}\right)^2-W_2\left(1-\frac{I^*}{I}\right)^2-W_3\left(1-\frac{R^*}{R}\right)^2$$
$$-W_4\left(1-\frac{B^*}{B}\right)^2+F(S,I,R,B)$$

where the function $F(S,I,R,B)$ is non-positive, now following the procedures by [20] and [21], we have;

$F(S,I,R,B)\leq 0$ for all $S,I,R,B$, Then $\frac{\mathrm{d}V}{\mathrm{d}t}\leq 0$ for all $S,I,R,B$ and it is zero when $S=S^*$, $I=I^*$, $R=R^*$, $B=B^*$ Hence the largest compact invariant set in $S,I,R,B$ such that $\frac{\mathrm{d}V}{\mathrm{d}t}=0$ is the singleton $E^*$ which is Endemic Equilibrium point of the model system (1).

Using LaSalles's invariant principle postulated by [22] we assert that $E^*$ is globally asymptotically stable in the interior of the region of $S,I,R,B$ and thus leads to the Theorem 6.3.

**Theorem 6.3** *If* $R_0>1$ *then the Typhoid Fever model system* (1) *has a unique endemic equilibrium point* $E^*$ *which is globally asymptotically stable in* $S,I,R,B$.

## 7. Numerical Analysis and Simulation

The section below presents the numerical analysis and simulation of the model, showing the behavior of Typhoid over a specific time period. **Figure 2** illustrates the human population dynamics during the disease outbreak. Clearly, early awareness and preventive measures can greatly reduce the impact of diseases like Typhoid. Information plays a key role in controlling diseases spread by human behavior and practices. The figure shows a rise in infections in the initial weeks, followed by a decline to an endemic equilibrium. This drop occurs as communities become aware of the disease and take precautions to limit its spread.

In the early weeks of an outbreak, the susceptible population drops exponentially to its endemic level due to a high infection rate. As infectious cases decline, recoveries also decrease proportionally. This trend is illustrated in **Figure 2**, showing a significant drop in $R$ as infection rates fall. Communities continue to live with diseases like typhoid because the bacteria's survival depends on environmental conditions. When conditions are unfavorable, the bacterial population declines exponentially but can rebound when conditions improve. **Figure 3** supports this, showing typhoid bacteria decreasing to its endemic equilibrium.

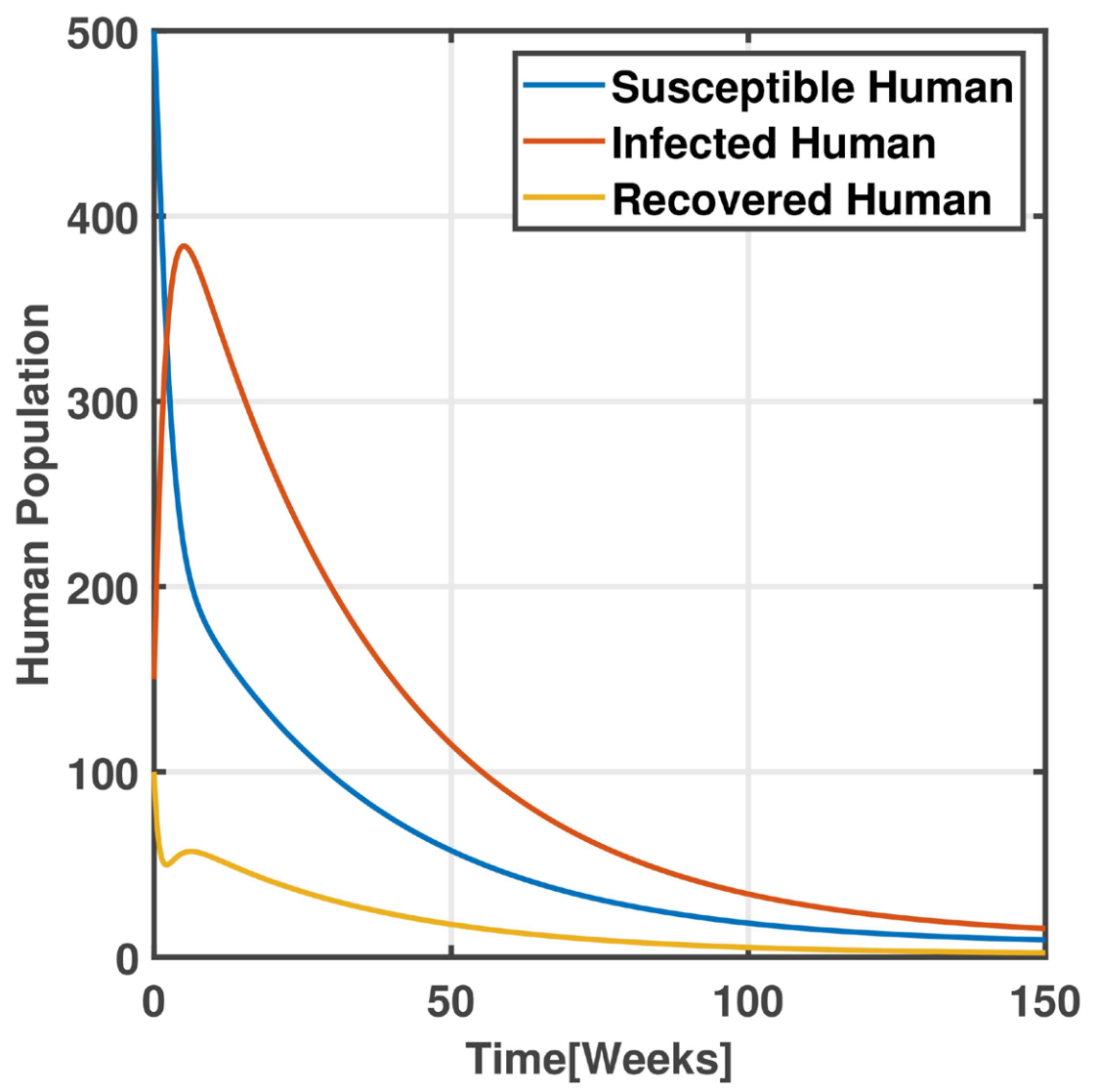

**Figure 2.** Dynamics of typhoid fever in human population.

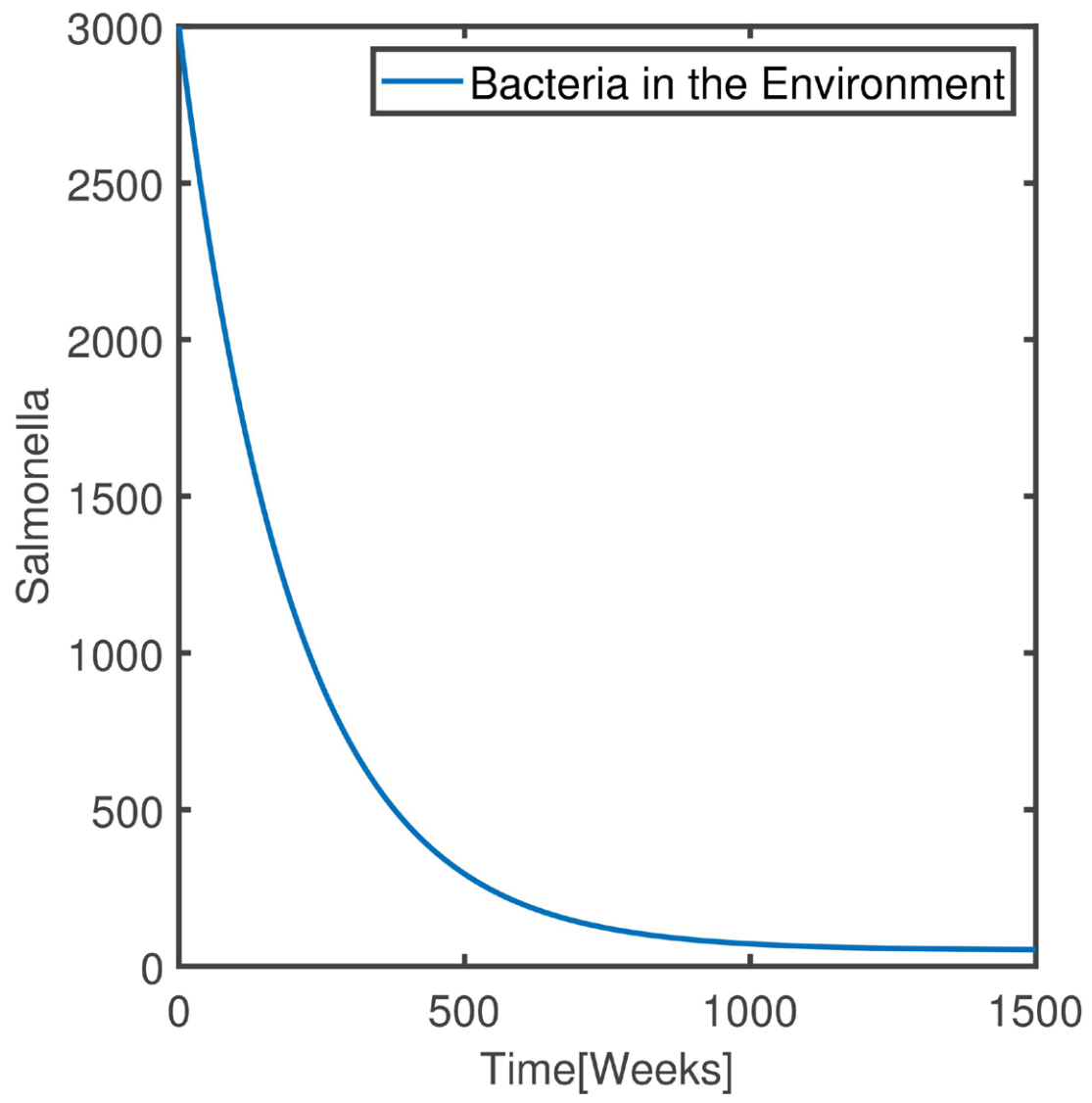

**Figure 3.** Dynamics of typhoid causing bacteria.

Typhoid fever can spread through direct contact between an infected and a susceptible person, but the main transmission route is through contaminated food and water. When environmental conditions support the growth of typhoid bacteria, disease prevalence tends to rise. Figure 4 illustrates that as environmental bacteria increase, the number of infectious humans also rises proportionally until it reaches a saturation point.

Moreover, Figure 4 points out the dependence on the number of bacteria in the environment and infectious human beings. As stated in the introduction above, infectious human beings shed typhoid bacteria in the environment. Thus, these two groups experience a mutual dependence in which the increase or decrease of one leads to the increase or decrease of the other.

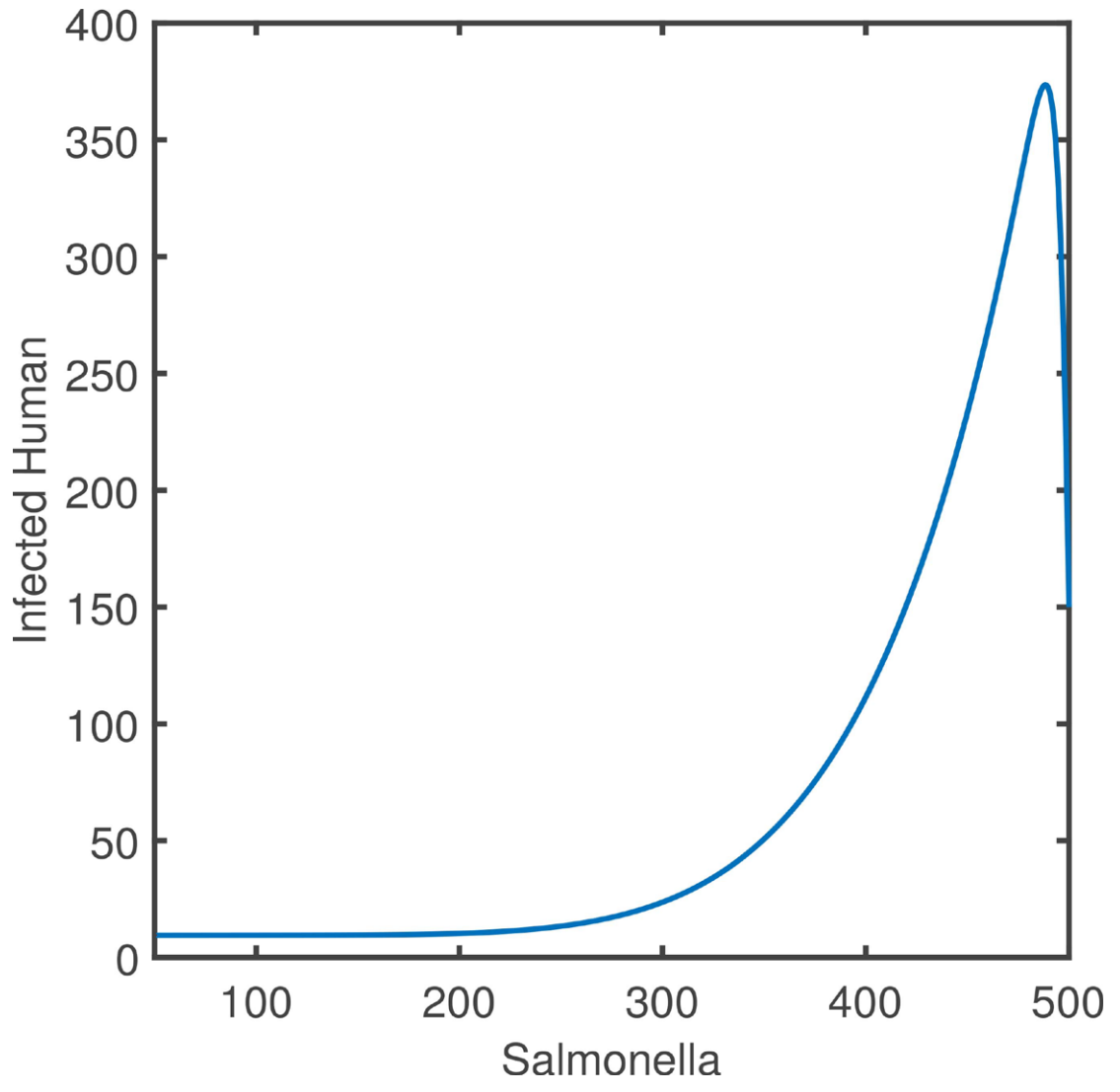

**Figure 4.** Dynamics of I with B.

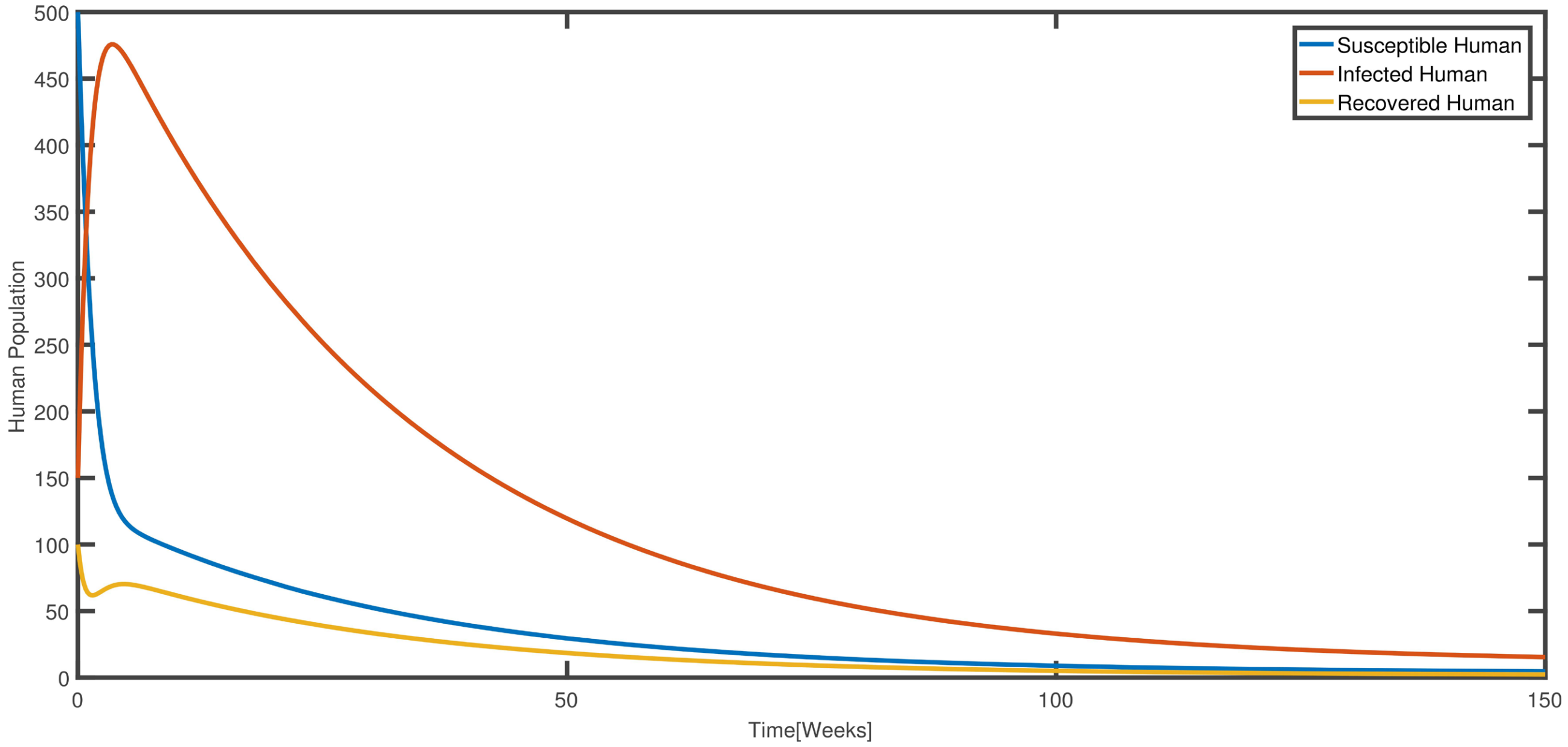

**Figure 5.** Role of information in typhoid transmission.

When the community is informed about risk behaviors that increase typhoid transmission, the spread of the disease is reduced. Figure 5 illustrates typhoid dynamics without such information or education. It shows a significant rise in infectious individuals compared to Figure 2, where information is provided. The results also indicate a longer infectious period without awareness, leading to greater disease impact and a higher risk of typhoid-related deaths.

## 8. Conclusion

This paper presents a model analyzing the role of information and education in the spread of typhoid fever. Conditions for both local and global stability are es-tablished. The basic reproduction number, derived using the next-generation matrix, quantifies the expected number of new infections caused by one infectious individual. Key factors influencing typhoid transmission include risk behavior, environmental bacterial concentration, and contact rate. The findings highlight the critical role of information and education in reducing risky behavior and emphasize their inclusion in effective typhoid control strategies.

## Contribution of the Study

The key contribution of this study is the formulation and analysis of a novel typhoid fever transmission model that incorporates information-based intervention as a behavioral factor influencing public response to disease outbreaks. The model captures how awareness campaigns and the dissemination of accurate health information can significantly reduce infection rates by encouraging preventive behavior. Additionally, we apply optimal control theory using Pontryagin Maximum Principle to determine the most effective strategies for minimizing the number of infected individuals and the cost of implementing interventions. Through theoretical analysis and numerical simulations, the study offers valuable insights into how information can be leveraged as a public health tool in combating typhoid fever.

## Acknowledgements

We thank everyone who supported this work directly or indirectly.

## Conflicts of Interest

The authors declare no conflicts of interest.